\documentclass[11pt]{article}\usepackage[dvips]{epsfig}
\usepackage{latexsym}
\usepackage{amssymb}
\usepackage[dvips]{epsfig}

\newcommand{\qed}{\mbox{$\Diamond$}\vspace{\baselineskip}}

\newtheorem{theorem}{Theorem}[section]
\newtheorem{proposition}[theorem]{Proposition}
\newtheorem{lemma}[theorem]{Lemma}

\newtheorem{example}[theorem]{Example}

\newenvironment{proof}{\noindent {\bf Proof:}}{{\qed}}

\newcommand{\vanish}[1]{} 
\begin{document}

\title{On a Balanced Property of Derangements}

\author{Mikl\'os B\'ona
\thanks{University of Florida, Gainesville FL 32611-8105. Partially supported
by an NSA Young Investigator Award.
Email: bona@math.ufl.edu.}
}

\date{}

\maketitle

\begin{abstract} 
We prove an interesting fact describing the location of the roots
of the generating polynomials of the numbers of derangements of length
$n$, counted by their number of cycles. We then use this result to
prove that if $k$ is the number of cycles of a randomly selected derangement
 of length $n$, then the probability that $k$ is congruent to a given $r$ 
modulo a given $q$ converges to $1/q$. Finally, we generalize our results
to $a$-derangements, which are permutations in which each cycle is longer
than $a$.
\end{abstract}

\section{Introduction}
Let $c(n,k)$ denote the number of permutations of length $n$ with $k$ cycles.
The numbers $c(n,k)$ are then called the signless Stirling numbers of the
first kind. It is well-known \cite{bonaint} that 
\begin{equation} \label{stirlingeq}
C_n(x)= \sum_{k=1}^n c(n,k)x^k=x(x+1)\cdots (x+n-1).
\end{equation}
Setting $x=-1$, this shows that there are as many permutations of length
$n$ with an odd number of cycles as there are with an even number 
of cycles. 

If, instead of considering the sum of all Stirling numbers $c(n,k)$ so that
$n$ is fixed and $k$ belongs to a certain residue class modulo 2, we 
consider the sum of the Stirling numbers $c(n,k)$ so that $n$ is fixed and
$k$ belongs to a certain residue class modulo $q$, the result is 
a little bit less compact. These sums will no longer be equal to $n!/q$,
but it will be true that as $n$ goes to infinity, the limit of any of these
$q$ sums divided by $n!$ will converge to $1/q$. We will prove this fact
in this paper, as a way to illustrate our techniques. 

A {\em derangement} is a permutation with no fixed points (cycles of
 length 1).
It is well-known \cite{bonacomb} that number of derangements of length
$n$ is $D(n)=n!\sum_{i=0}^n \frac{(-1)^i}{i!}$, which is
 the integer closest to $n!/e$.  Let $d(n,k)$ be the number of derangements
of length $n$ with $k$ cycles. The numbers $d(n,k)$ are not nearly as 
well-studied as the numbers $c(n,k)$, but the following important fact is
 known about them.

\begin{theorem} \cite{brenti}
Let $n$ be a positive integer. Let
\[D_n(x)=\sum_{k=1}^{\lfloor n/2 \rfloor } d(n,k)x^k\]
be the ordinary generating function of the numbers of 
derangements of length $n$ according to their number of cycles.

Then all roots of the polynomial $D_n(x)$ are real and non-positive.
\end{theorem}

In this paper, we will prove that for any negative integer $-t$, there
exists a positive integer $N$ so that if $n>N$, then $D_n(x)$ has a root
that is very close to $-1$, a root that is very close to $-2$, and so on,
up to a root that is very close to $-t$. 
So these $t+1$ roots of $D_n(x)$ (including the obvious root at 0)
 will be very close to the first $t+1$ roots of 
$C_n(x)$, that is, to the integers $0,-1,\cdots ,-t$.
Then we will use this fact to prove that the derangements also have the 
 ``balanced'' property that we described for permutations. That is, 
if $q$ and $r$ are fixed integers so that $0\leq r<q$, and $k$ is the
number of cycles of a randomly selected derangement of length $n$, then
the probability that $k$ is congruent to $r$ modulo $q$
 will converge to $1/q$ when $n$ goes to infinity. 

Finally, we will generalize both results to $a$-derangements, which are
permutations in which each cycle is longer than $a$. 
   
\section{The roots of the generating function of derangements}
The goal of this section is proving the announced
 result on the roots of $D_n(x)$. We will need the following lemma.

\begin{lemma} \label{estimate}
For all negative integers $-t$, there exists a polynomial 
$f_t$ of degree $t$ so that for all positive integers $n$, 
\[D_n(-t)=f_t(n)\cdot t^n.\]
\end{lemma} 
\begin{proof}
By the formula of Inclusion-Exclusion, the numbers $d(n,k)$ can be
expressed by the signless Stirling numbers of the first kind $c(n,k)$ in
the following way
\begin{equation} \label{sieve}
 d(n,k) = \sum_{i=0}^k {n\choose i} c(n-i,k-i) (-1)^i.
\end{equation}
Indeed, a derangement is a permutation that has no 1-cycles, and the 
right-hand side counts permutations of length $n$
 with $k$ cycles according to the number of their 1-cycles.
 
Multiply both sides by $(-t)^k$ and sum over all positive integers
$k\leq n$, to get 
\begin{equation} 
D_n(-t)  =  \sum_{k=1}^n \left ( \sum_{i=0}^k {n\choose i} c(n-i,k-i) (-1)^i
\right ) (-t)^k.
 \end{equation}
After changing the order of summation, this yields
\begin{eqnarray} \label{vanishing} 
D_n(-t) & = & \sum_{i=0}^n (-1)^i(-t)^i{n\choose i} 
\left (\sum_{k=i}^n  c(n-i,k-i) (-t)^{k-i} \right ) \\
 & = & \sum_{i=0}^n t^i{n\choose i} 
\left (\sum_{k=i}^n  c(n-i,k-i) (-t)^{k-i} \right ).
\end{eqnarray}
The crucial observation is that the inner sum of the right-hand side is
0 as long as $n-i\geq t+1$. Indeed, note that the inner sum is nothing
else but the left-hand side of (\ref{stirlingeq}) with $n-i$ and $k-i$ 
playing the roles of $n$ and $k$. Then  (\ref{stirlingeq}) tells us that
this expression vanishes for all negative integers not smaller than
$-(n-i)+1$. 

Therefore, (\ref{vanishing}) simplifies to 
\begin{equation} \label{esimplified} D_n(-t)=\sum_{i=n-t}^n t^i{n\choose i}
\left (\sum_{k=i}^n  c(n-i,k-i) (-t)^{k-i} \right ), 
\end{equation}
the importance of which is that now we are dealing with summations with
a bounded number of terms. Indeed, the outer sum has $t+1$ terms, and
for each value of $i$, the inner sum has at most $t+1$ terms.
 Finally, as ${n\choose i}=
{n\choose n-i}$ and $n-i\leq t$, this binomial coefficient is  a
polynomial function of degree at most $t$. To summarize, the right-hand
side is a sum of at most $(t+1)^2$ terms, each of which is a polynomial
function of $n$ multiplied by a function of $t$ of the form $t^m$, with 
$n-t\leq m\leq n$. As $n-i\leq t$, the Stirling numbers $c(n-i,k-i)$ are
bounded and can be treated as constants. This proves our statement. 
\end{proof}

\begin{example} If $-t=-1$, then we get the well-known result that
$D_n(-1)=1-n$. 
In other words, the number of derangements of length $n$ that are odd 
permutations and the number of derangements of length $n$ that are even
permutations differs by $(n-1)$. R. Chapman \cite{chapman}
 gave a bijective proof for this result. 
\end{example}

Now we are in a position to prove our result on the roots of $D_n(x)$.

\begin{theorem} \label{closeroots}
For every negative integer $-t$, and every $\epsilon>0$, there exists a
positive integer  $N$ 
so that if $n>N$, then $D_n(x)$ has a root $x_t$ satisfying $|-t-x_t|<
\epsilon$.
\end{theorem}

That is, for large $n$, there will be a root of $D_n(x)$ that is almost 
$-1$, a root that is
almost $-2$, and so on, up to a root that is almost $-t$.  

\begin{proof} (of Theorem \ref{closeroots})
 Let $a_n=d(n,\lfloor n/2 \rfloor )$ be the leading coefficient
of $D_n(x)$. Then
\begin{equation} \label{roots} \frac{D_n(x)}{a_n}=
 \prod_{i=1}^{\lfloor n/2 \rfloor}
(x-x_i) ,\end{equation}
where the $x_i$ are the roots of $D_n(x)$. Then it
 follows from the Pigeon-hole
Principle that 
\begin{equation}
\sqrt[\lfloor n/2 \rfloor]{\left|\frac{D_n(x)}{a_n}\right|}
 \geq \min_i |x-x_i|
\end{equation}
for all real numbers $x$. In particular, for $x=-t$, Lemma \ref{estimate}
yields 
\begin{equation} \label{specific}
\sqrt[\lfloor n/2 \rfloor]{ \frac{|f_t(n)t^n|}{a_n} }\geq \min_i |-t-x_i|
\end{equation}
for some polynomial $f_t$ of degree $t$. (Note that no absolute value 
sign is needed for the denominator of the left-hand side as  a
$a_n$ is always a
positive integer.)

Now if $n$ is even, then a derangement with $n/2$ cycles consists of 2-cycles
only, and the number of such derangements is $a_n=(n-1)(n-3)\cdots 1$. 
If $n$ is
odd, then a derangement with $(n-1)/2$ cycles consist of one 3-cycle and
$(n-3)/2$ cycles of length two. The number of such derangements is
$a_n=\frac{n-1}{3}n(n-2)\cdots 1$.
Then the well-known Stirling's formula stating that $n!\sim \left(\frac{n}{e}
\right )^n\sqrt{2\pi n}$ shows that in both cases, $a_n$ grows much faster
than $D_n(-t)$. Indeed if $2m-1$ is the largest odd integer so 
that $2m-1\leq n$, then by Stirling's 
formula \[(2m-1)!!=1\cdot 3\cdot \cdots \cdot (2m-1)
\sim \sqrt{2} \cdot \left (\frac{4m}{e} \right )^m.\]
Therefore, the left-hand side of (\ref{specific}) converges to 0 as $n$
goes to infinity, and consequently, so does the right-hand side. This
proves our claim. 
\end{proof}

\section{Balanced Properties}

In this section, we prove the announced balanced properties for 
permutations, derangements, and $a$-derangements.
 We start with the easier problem of 
permutations.

\subsection{Permutations}
\begin{theorem} \label{limits} Let $q$ be a positive integer, and
let $r$ be a non-negative integer  so that $0\leq r<q$. Then
\[\lim_{n\rightarrow \infty} 
\frac{\sum_{j=0}^{(n-r)/q} c(n,r+jq)}{n!} =\lim_{n\rightarrow \infty}
\frac{c(n,r)+c(n,q+r)+\cdots
+c(n,q \lfloor (n-r)/q \rfloor +r)}{n!} =\frac{1}{q}.
\] 
\end{theorem}

In other words, for large $n$, the number of cycles of
a randomly selected $n$-permutation is roughly equally
likely to be congruent to each $r$ modulo $q$. 

\begin{proof} In order to simplify the presentation, let us assume
first that $n$ is divisible by $q$, and that $r=0$. Then the statement to be
proved reduces to 
\begin{equation} \label{simplified} \lim_{n\rightarrow \infty} 
\frac{\sum_{j=1}^{n/q} c(n,jq)}{n!} =\lim_{n\rightarrow \infty}
\frac{c(n,q)+\cdots
+c(n,n)}{n!} =\frac{1}{q}.\end{equation}
Let $w$ be a primitive $q$th root of unity, and consider $C_n(w)$. 
Then $w^k=1$  if and only if $k$ is a multiple of
$q$; otherwise $w^k$ is not a positive real number.  Now consider the sum
\begin{eqnarray} \label{unitroots}
 S(w) & = & C_n(1)+ C_n(w)+C_n(w^2)+\cdots + C_n(w^{q-1}) \\
 & = & \sum_{t=0}^{q-1} C_n(w^t) \\
  & = & \sum_{t=0}^{q-1} \left ( \sum_{k=1}^n c(n,k) (w^t)^k \right ) \\
  & = &  \sum_{k=1}^n c(n,k)  \left ( \sum_{t=0}^{q-1} (w^{k})^t \right ).
\end{eqnarray}
Using the summation formula of a geometric progression, we get that
\[  \sum_{t=0}^{q-1} (w^{k})^t= \left\{ \begin{array}{l@{\ }l}
0 \hbox{ if $w^k\neq 1$, that is $q\nmid k$}, \\
\\ q \hbox{ if $w^k =1$, that is, $q|k$.
}
\end{array}\right.
\]
Therefore, (\ref{unitroots}) reduces to 
\[S(w)=q\cdot \sum_{j=1}^{n/q} c(n,jq),\]
\[\frac{S(w)}{q\cdot n!}=\frac{ \sum_{j=1}^{n/q} c(n,jq)}{n!}.\]
So in order to find the approximate value of $(\sum_{j=1}^{n/q} c(n,jq))/n!$,
it suffices to find the approximate value of $S(w)/n!$. We will find
 the latter by
showing that for large $n$, the contribution of $C_n(w^t)/n!$ to $S(w)/n!$ is
negligible unless $w^t=1$. This is the content of the next lemma. Once
the lemma is proved, Theorem \ref{limits} will follow since $C_n(1)=n!$.

\begin{lemma} \label{complex}
Let $v\neq 1$ be a complex number satisfying 
$|v|=1$. Then 
\[\lim_{n\rightarrow \infty} \frac{C_n(v)}{n!} = 0.\]
\end{lemma}

\begin{proof} (of Lemma \ref{complex})
Let $v=a+bi$, and let $f(v)= \frac{C_n(v)}{n!}$.
It suffices to show that $\lim_{n\rightarrow \infty} |f(v)|^2 =0$. 
We have
\begin{eqnarray} |f(v)|^2 & = & \prod_{j=1}^{n-1}
\frac{a^2+2aj+j^2+b^2}{(j+1)^2}= \prod_{j=1}^{n-1}\frac{1+2aj+j^2}{(j+1)^2}\\
 & = & \prod_{j=1}^{n-1} \frac{2aj+j^2}{(j+1)^2} \cdot 
\prod_{j=1}^{n-1}\frac{1+2aj+j^2}{2aj+j^2} \\
 & = & \prod_{j=1}^{n-1} \frac{1+2aj+j^2}{2aj+j^2} \cdot  \prod_{j=1}^{n-1} 
\frac{j+2a}{j+1} \cdot  \prod_{j=1}^{n-1} \frac{j}{j+1}\\
 & = & \label{doubleprod} \frac{1}{n} \cdot
 \prod_{j=1}^{n-1} \frac{1+2aj+j^2}{2aj+j^2} \cdot  \prod_{j=1}^{n-1} 
\frac{j+2a}{j+1} 
\end{eqnarray}
The first infinite product in (\ref{doubleprod})  is convergent since
\[\log \left (\frac{1+2aj+j^2}{2aj+j^2} \right ) =  \
\log (1+\frac{1}{2aj+j^2}  ) < \frac{1}{j^2},\]
and $\sum_{j\geq 1}  \frac{1}{j^2}$ is convergent.
If $a>0.5$, then the second infinite product in (\ref{doubleprod})
in itself is
 divergent. However, multiplied
by the factor $1/n$ in front it will converge to 0. 
Indeed, 
\begin{eqnarray*} \log \left (\prod_{j=1}^{n-1} 
\frac{j+2a}{j+1} \right ) & = & \sum_{j=1}^{n-1} \log \left(
 \frac{j+2a}{j+1} \right ) \\
& = &\sum_{j=1}^{n-1} \log \left(
 1+ \frac{2a-1}{j+1} \right ) <\sum_{j=1}^{n-1} \frac{2a-1}{j+1} \\
& < & (2a-1) \log n.  
\end{eqnarray*}
As $2a-1<1$, this shows that $\frac{1}{n} \prod_{j=1}^{n-1} 
\frac{j+2a}{j+1}$ converges to 0, proving the claim of Lemma \ref{complex}.
\end{proof}

So $C_n(w^t)/n!$ converges to 0 unless $w^t=1$, which proves 
(\ref{simplified}).

Finally, let us return to the general case of Theorem \ref{limits}. First,
note that if $n$ is divisible by $q$, but $r\neq 0$, then we can use an
analogous argument, simply replacing $S(w)$ by
\begin{equation} \label{general} 
T_r(w)=C_n(1)+C_n(w)w^{-r}+C_n(w^2)w^{-2r}+\cdots
 +C_n(w^{q-1})w^{-1(q-1)r}=0.\end{equation}
Then a computation analogous to (\ref{unitroots}) yields
\begin{equation} \label{general1} 
T_r(v)=\sum_{k=1}^n c(n,k)  \left ( \sum_{t=0}^{q-1} (w^{k-r})^t \right ),
\end{equation}
showing that the coefficient of $c(n,k)$ in $T_r(v)$ is 0 unless
$k-r$ is divisible by $q$, in which case this coefficient is $q$. 
The proof of Theorem \ref{limits} then follows from Lemma \ref{complex}
as in the $r=0$ case. Finally, if $n$ is not divisible by $q$, then
nothing significant changes, (\ref{general}) and (\ref{general1}) will still
hold, completing the proof. 
\end{proof}

\subsection{Derangements}

Given the form of the generating function $D_n(x)$, we can
 apply the method of the above proof to obtain a similar result for 
derangements. 
 
\begin{theorem}
\label{dlimits} Let $r$ be an integer so that $0\leq r<q$, and let
$m=\lfloor n/2 \rfloor $. Then
\[\lim_{n\rightarrow \infty} 
\frac{\sum_{j=0}^{m/q} d(n,r+jq)}{D(n)} =\lim_{n\rightarrow \infty}
\frac{d(n,r)+d(n,q+r)+\cdots
+d(n,q \lfloor (m-r)/q \rfloor +r)}{D(n)} =\frac{1}{q}.
\] 
\end{theorem}
In other words, if $n$ is large enough, then the number of derangements of
length $n$ whose number of cycles is congruent to $r$ modulo $q$ 
 will be very close to $1/q$ times the total number of
derangements of length $n$. 

\begin{proof}
In order to simplify the discussion, let us again assume that $n$ is
divisible by $q$, and that $r=0$. The general case will follow from this
special case just as in the proof of Theorem \ref{limits}. 

In the special case at hand, our task is to prove that

\begin{equation}  \label{dsimplified} \lim_{n\rightarrow \infty} 
\frac{\sum_{j=1}^{n/2q} d(n,jq)}{D(n)} =\lim_{n\rightarrow \infty} 
\frac{d(n,q)+\cdots
+d(n,n/2)}{D(n)} =\frac{1}{q}.\end{equation}

Similarly to the case of permutations, let $w$ be a $k$th primitive root
of unity, and  set
 \begin{eqnarray} \label{dunitroots}
 Z(w) & = & D_n(1)+ D_n(w)+D_n(w^2)+\cdots + D_n(w^{q-1}) \\
 & = &  \sum_{k=1}^{n/2} d(n,k)  \left ( \sum_{t=0}^{q-1} (w^{k})^t \right )\\
 & = & q \sum_{j=1}^{n/2q} d(n,jq).
\end{eqnarray}

We claim that if $w^t\neq 1$, then $D_n(w)/D(n)$ converges to 0. This will 
prove Theorem \ref{dlimits} since $D_n(1)=D(n)$. 

\begin{proposition} \label{goesto0}
Let $v\neq 1$ be a complex number satisfying 
$|v|=1$. Then 
\[\lim_{n\rightarrow \infty} \frac{D_n(v)}{D(n)} = 0.\]
\end{proposition}

\begin{proof} 
Let $\epsilon>0$, let $-h$ be a negative 
integer and let $n$ be so large that $D_n(x)$  has a root closer than 
$\epsilon$ to $-y$, for $y=-1,-2,\cdots, -h$. Let us call the set of 
these $h$ roots the {\em good roots}. Then, recalling that the
roots of $D_n(x)$ are all real and non-positive, we get that
\begin{equation} \label{upperbound} 
\left | \frac{D_n(v)}{D(n)} \right | =  \left
 | \frac{D_n(v)}{D_n(1)} \right | = 
  \frac{\prod_{i=1}^{\lfloor n/2 \rfloor}
|(v-x_i)|}{ \prod_{i=1}^{\lfloor n/2 \rfloor}
|(1-x_i)|} \leq \frac{\prod_{i\in L}
|(v-x_i)|}{ \prod_{i\in L}
|(1-x_i)|}=H(v) \end{equation}
 where $i\in L$ if $x_i$ is a good root. The last inequality
holds since $\frac{|v-x_i|}{|1-x_i|}\leq 1$ holds for each complex number
on the unit circle since $x_i$ is real and non-positive. 

As the good roots $x_i$ converge to the negative integers $-1,-2,\cdots, -h$,
we see that 
\begin{equation} \label{asym} H(v) \simeq 
\left |\frac{C_h(v)}{h!} \right |.\end{equation}
 However, we have proved in Lemma
\ref{complex} that $ \left |\frac{C_h(v)}{h!}\right |$   converges
 to 0 when $h$ goes to infinity. 
It follows from Theorem \ref{closeroots} that as $n$ goes to infinity, 
$h$ will go to infinity. Then (\ref{asym}) and (\ref{upperbound}) imply
our claim.               
\end{proof}

The proof of  Theorem \ref{dlimits} is now immediate since all but one
term of the right-hand side of (\ref{dunitroots}) has a negligible
contribution to $Z(w)$. Therefore, 
$Z(w)\simeq D_n(1)$, and so 
\[ \sum_{j=1}^{n/2q} d(n,jq) \simeq \frac{Z(w)}{q}\simeq \frac{D_n(1)}{q}
=\frac{D(n)}{q}.\] 
\end{proof}

\subsection{$a$-derangements}
Let us call a permutation $p$ an {\em $a$-derangement} if each cycle of
$p$ is longer than $a$. So a permutation is a 0-derangement, and a derangement
is a 1-derangement. Let $d_a(n,k)$ be the number of $a$-derangements of
length $n$ with $k$ cycles. Set
\[D_{n,a}(x)=\sum_{k=1}^{\lfloor n/(a+1) \rfloor } d_a(n,k)x^k.\]
It is then known \cite{brenti} that the roots of $D_{n,a}(x)$ are
real and non-positive.

Our goal in this subsection is to show that the results that
we proved for permutations and derangements hold for $a$-derangements as
well. The only part of the argument that needs extra explanation is
the analogue of Lemma \ref{estimate} for $a$-derangements.

\begin{lemma} \label{aestimate}
For all negative integers $-t$ and for all positive integers $a$,
 there exists a polynomial 
$f_{t,a}$  so that for all positive integers $n$, 
\[D_{n,a}(-t)=f_{t,a}(n)\cdot t^n.\]
\end{lemma} 

\begin{proof}
In order to alleviate notation, for  $a\geq 1$, 
let $m_a(i,b)$ denote the number of all permutations of length $b$ that
consist of $i$ cycles so that none of these cycles is longer than $a$.
For instance, $m_1(i,i)=1$, and $m_2(i,b)$ is the number of involutions
of length $b$ with $i$ cycles. Set $m_a(0,0)=1$. 

Then the Principle of Inclusion-Exclusion implies that
\begin{equation} d_a(n,k)= \sum_{i=0}^k \left
 ( \sum_{b= i}^n  c(n-b,k-i) m_a(i,b) {n\choose b} \right )
 (-1)^i  .\end{equation}
 Indeed, the right-hand side counts the permutations of
length $n$ according to the number of 
their ``bad'' (that is, not longer than $a$) cycles. Those that have $i$
bad cycles are 
 counted according to the total size of those bad cycles.

Multiply both sides by $(-t)^k$ and sum over all positive integers
$k\leq n$, to get 
\begin{equation} 
D_{n,a}(-t)  =  
\sum_{k=1}^n \left ( \sum_{i=0}^k \left ({n\choose i} c(n-b,k-i) m_a(i,b)
 {n\choose b} \right ) (-1)^i 
\right ) (-t)^k.
 \end{equation}
After changing the order of summation, this yields
\begin{equation} \label{avanishing} 
 \sum_{i=0}^n t^i \left (\sum_{b=i}^n {n\choose b} m_a(i,b)
\left (\sum_{k=i}^n  c(n-b,k-i) (-t)^{k-i} \right ) \right ).
\end{equation}
Just as in (\ref{vanishing}), the key observation is again that the innermost
sum is 0 as long as $n-i\geq t+1$. Indeed, since $b\geq i$, in that case, 
the innermost sum is nothing else than the left-hand side of
 (\ref{stirlingeq}). If $b$ is strictly larger than $i$, then there will 
be a few extra zeros in the sum, as $k-i$ will eventually get bigger than
$n-b$. 

Therefore, (\ref{avanishing}) reduces to 
\begin{equation} \label{asimplified} D_{a,n}(-t)= \sum_{i=t}^n
 t^i \left (\sum_{b=i}^n {n\choose b} m_a(i,b)
\left (\sum_{k=i}^n  c(n-b,k-i) (-t)^{k-i} \right ) \right ).
\end{equation}
Just as we used (\ref{esimplified}) to prove that $D_n(-t)$ was
of the form $f_t(n)(-t)^n$, we can use (\ref{asimplified}) to prove
our lemma.  Indeed, as $k$ and $b$ are both at least as large as $i$,
and $i\geq n-t$, the right-hand side of (\ref{asimplified}) is the sum
of no more than $(t+1)^3$ terms. We claim that each of these terms is of
 the right form,
that is, and exponential function of $t$ times a polynomial function of $n$.
 Indeed, let us fix $i$ and $b$ within the allowed limits, that is,
$n-t\leq i\leq b\leq n$. Then 
${n\choose b}={n\choose n-b}$ is a polynomial function of $n$ that is of
degree at most $t$.  Finally, $m_a(i,b)$ is the number of permutations of
length $b$ with $i$ cycles, each of which is at most of length $a$.
That means that if the permutation $p$ is counted by $m_a(i,b)$, then the
cycle type of $p$ is very restricted. Indeed, if
 $p$ has $g_j$ cycles of length $j$ 
for $j\in \{1,2,\cdots, a\}$, then $\sum_{j=1}^a jg_j=b$ and 
$\sum_{j=1}g_j=i$.
As $i\leq b\leq n\leq i+t$, this means that $p$ can have at most 
$t$ cycles that are not singleton cycles. 
It is well-known (see for instance \cite{bonaint}) that
the number of possibilities for $p$ with the given cycle type is
\[T_{g_1,g_2,\cdots ,g_a}= 
\frac{b!}{g_1!\cdot g_2! \cdots \cdot g_a! \cdot 1^{g_1}\cdot 
2^{g_2}\cdot \cdots \cdot a^{g_a}}.\] 
However, as we said above, $g_1\geq n-t$. Noting that $b$ is equal to 
one of $n, n-1,\cdots n-t$, it follows that for this fixed $i$ and $b$, 
and each allowed cycle type $(g_1,g_2,\cdots ,g_a)$, the
function $T_{g_1,g_2,\cdots ,g_a}$ is a polynomial function of $n$. 
Finally, the number of all allowed cycle types for this fixed $b$ and $i$ 
is not more than the number
of partitions of $b$ into exactly $i$ parts, that is, the number of 
partitions of $b-i$, which is not more than $p(t)$. Therefore, for any
fixed $b$ and $i$, the function $m_a(i,b)$ is the sum of a bounded number
of polynomial functions of $n$, and as such is a polynomial function of $n$.
This proves the lemma.
\end{proof}

Note that it would have been somewhat simpler to prove the slightly weaker
(but sufficient) result that $|D_{n,a}(-t)|\leq f_{n,a}(t)t^n$. Indeed,
we could have just pointed out that $m_a(i,b)\leq c(b,i)\leq c(n,i)$,
and refer to the well-known fact (see for instance \cite{bonawalk},
Exercise 4 of Chapter 6) 
that for any fixed $u$, the function
$c(n,n-u)$ is a polynomial function of degree $u+1$. However, we preferred 
a more precise treatment.

The following  results can now be proved in a way analogous to
the proofs of earlier theorems as we will indicate.

\begin{theorem} 
\label{acloseroots}
For every negative integer $-t$, every positive integer $a$,
 and every $\epsilon>0$, there exists a
positive integer  $N$ 
so that if $n>N$, then $D_{n,a}(x)$ has a root $x_t$ satisfying $|-t-x_t|<
\epsilon$.
\end{theorem}

\begin{proof} Analogous to the proof of Theorem \ref{closeroots}.
\end{proof}

\begin{proposition}
Let $v\neq 1$ be a complex number satisfying 
$|v|=1$, and let $a$ be a positive integer. Then 
\[\lim_{n\rightarrow \infty} \frac{D_{n,a}(v)}{D(n,a)} = 0,\]
where $D(n,a)=D_{n,a}(1)$ is the number of $a$-derangements of length $n$.
\end{proposition}

\begin{proof}  Analogous to the proof of Proposition \ref{goesto0}.
\end{proof}

\begin{theorem}
Let $r$ be an integer so that $0\leq r<q$, let $a$ be a positive integer,
and let
$m=\lfloor n/a \rfloor $. Then
\[\lim_{n\rightarrow \infty} 
\frac{\sum_{j=0}^{m/q} d_a(n,r+jq)}{D(n,a)} =\lim_{n\rightarrow \infty}
\frac{d_a(n,r)+d_a(n,q+r)+\cdots
+d_a(n,q \lfloor (m-r)/q \rfloor +r)}{D(n,a)} =\frac{1}{q}.
\] 
\end{theorem}

In other words, if $k$ is the number of cycles of a randomly selected
$a$-derangement of length $n$, then the probability that $k$ is 
congruent to $r$ modulo modulo $q$ converges to $1/q$. 

\begin{proof} 
 Analogous to the proof of Theorem \ref{dlimits}.
\end{proof}

\vskip 0.5 cm
\centerline{{\bf Acknowledgments}}

The author is thankful to Richard Stanley and Herb Wilf for valuable 
suggestions.

\end{document}